\title{Topology of symplectomorphism groups and ball-swappings}

\documentclass[11pt]{article}

\usepackage{amsfonts}
\usepackage{amsthm}
\usepackage{latexsym, amssymb, bbm}
\usepackage{xcolor}
\usepackage{graphicx}
\usepackage{amsmath}
\usepackage{enumerate}
\usepackage[all]{xy}
\usepackage[margin=1.3in]{geometry}
\usepackage[mathscr]{eucal}
\usepackage{latexsym, amssymb, bbm, mathrsfs,mathabx}
 \usepackage{caption} \usepackage{subcaption}

\usepackage{amsmath,amsthm,amsfonts,amssymb,graphicx,psfrag}
\usepackage{color}

\usepackage{hyperref}

\usepackage{tikz}
\usepackage{comment}

\newtheorem{thm}{Theorem}[section]

\newtheorem{prop}[thm]{Proposition}
\newtheorem{lemma}[thm]{Lemma}
\newtheorem{rmk}[thm]{Remark}

\newtheorem{defn}[thm]{Definition}
\newtheorem{corr}[thm]{Corollary}

\newtheorem{conj}[thm]{Conjecture}

\theoremstyle{definition}
\newtheorem{example}[thm]{Example}

\DeclareMathOperator{\Diff}{Diff}

\newcommand{\wt}[1]{\widetilde{#1}}

\newcommand{\ov}[1]{\overline{#1}}

\newcommand{\bdf}{\begin{defn}}
\newcommand{\edf}{\end{defn}}

\newcommand{\bthm}{\begin{thm}}
\newcommand{\ethm}{\end{thm}}

\newcommand{\blem}{\begin{lemma}}
\newcommand{\elem}{\end{lemma}}

\newcommand{\bcor}{\begin{corr}}
\newcommand{\ecor}{\end{corr}}

\newcommand{\bprop}{\begin{prop}}
\newcommand{\eprop}{\end{prop}}

\newcommand{\brmk}{\begin{rmk}}
\newcommand{\ermk}{\end{rmk}}

\newcommand{\bpf}{\begin{proof}}
\newcommand{\epf}{\end{proof}}

\newcommand{\bex}{\begin{example}}
\newcommand{\eex}{\end{example}}

\numberwithin{equation}{section}

\def\eL{\EuScript{L}}

\def\C{\mathbb{C}}

\def\Q{\mathbb{Q}}

\def\Z{\mathbb{Z}}

\def\CP{\mathbb{CP}}
\def\RP{\mathbb{RP}}

\def\cC{\mathcal{C}}

\def\fuk{\EuScript{F}uk}
\def\ov{\overline}

\def\w{\omega}

\def\xkm2{\overline{X}_{k-2}}

\newtheorem{prob}[thm]{Problem}

\newcommand{\Conf}{Conf^{ord}}
\newcommand{\mJ}{\mathcal{J}}

\begin{document}

\author{Jun Li\thanks{J.L. is partially supported by NSF grant 
DMS-1611680. } and Weiwei Wu\thanks{W.W. is partially
supported by Simons Collaboration Grant 524427.}}

\AtEndDocument{\bigskip{\footnotesize
  \textsc{Jun Li, Department of Mathematics, University of Michigan, Ann Arbor, MI 48109} \par
  \textit{E-mail address}: \texttt{lijungeo@umich.edu} \par

}\bigskip{\footnotesize
  \textsc{Weiwei Wu, Department of Mathematics,
      University of Georgia,
      Athens, GA 30606} \par
  \textit{E-mail address}: \texttt{weiwei.wu@uga.edu} \par

}}

\date{\today}
\maketitle

\begin{abstract}

We summarize some recent progress and problems on the symplectomorphism groups, with an emphasis on the connection to the space of ball-packings.

\end{abstract}

\section{Introduction} 
\label{sec:introduction}

Let $X$ be a symplectic manifold, with symplectic form $\w$. A
symplectomorphism is a diffeomorphism $\phi: X \rightarrow X$
such that $\phi^*\w = \w$. The main theme of this article is the group of symplectomorphisms $Symp(M,\w)$, which can be endowed with a $C^{\infty}$-topology and becomes an infinite-dimensional Lie group.  Since Gromov's seminal work in \cite{Gro85}, the topology of symplectomorphism groups has been an active field of study.  As of now, the subject has developed into a large industry, and many different faces have been discovered and explored.

One of the most attractive problems on $Symp(X)$, after Gromov \cite{Gro85}, is to understand its (weak) homotopy type of $Symp(X)$.  In general, we know $Symp(X)$ is always a countable CW-complex \cite{MS98}, but any concrete computations are highly nontrivial, see for example \cite{Abr98,McD00,AM00,Anj02,LP04,Buse11}.  This thread of study of the full homotopy type has been deeply intertwined with the stratification of almost complex structures, which again inspired a series of interesting works \cite{McD00, AGK09, Buse05, LL16}, etc.  If one narrows the attention to the connectedness of $Symp_h(X)$, the symplectomorphism subgroup that acts homologically trivially on $X$, interesting phenomenon can already be observed.  Seidel first constructed the Lagrangian Dehn twists as symplectomorphisms, and proved that they are smoothly isotopic to identity but not symplectically isotopic to identity in dimension $4$  in many cases \cite{Se99,Sei08}.  He further pointed out the connection between such symplectomorphisms and the singularity theory which goes back to Arnold \cite{Se4dim, Seithesis}.  This thread of ideas leads to the symplectic Picard-Lefschetz theory, which shapes a large part of homological mirror symmetry, but is beyond the scope of our discussion here.

The perspective of the current article on $Symp(X)$ is an attempt to extend Seidel's works in dimension $4$ and make connections with yet another classical problem in symplectic geometry: the ball-packing problem.  Since the discovery of Gromov's celebrated non-squeezing theorem, the ball-packing problem has been one of the central themes in symplectic geometry, see \cite{Bir97,Bi97,Bi96,Bir01, Opsh07, MS12, CCDDHR} for some examples of the $4$ dimensional case.  The higher dimension cases remain one of the most exciting areas in current research, see \cite{Pel08, BH11, BH13, HK14, HK18} for just a few examples.  

To make the connection to symplectomorphism groups, we recall the symplectic description of blow-ups \cite{MP94}.  Take a symplectic $4$ manifold $X$ and a ball packing $\phi:\coprod_{i=1}^k B(c_i)\hookrightarrow X$, one may obtain a symplectic $k$-fold blowup of $X$ by quotienting the Hopf fibration of the boundary of the complement of the balls.  Seidel pointed out one should consider a ``universal'' fibration of $X\#k\ov\CP^2$ over the configuration space of points in $X$, whose description is recalled in Section \ref{BS}.  The monodromy of this fibration, governed by the fundamental group of the configuration space, contains a large amount of information of the \textit{symplectic mapping class group} $\pi_0(Symp(X\#k\ov\CP^2))$.

Another similar but more topological perspective is based on Mcduff's observation that the space of ball-packing is connected \cite{McD96}, preceded by the work of Biran \cite{Bi96}.  In special cases, this allows one to study the symplectomorphism group via the homotopy fibration \eqref{e:action}, as is explained in Section \ref{sub:symplectic_mapping_class_groups_of_rational_surfaces}.
The second author defines a kind of possible generators of $\pi_0Symp(M,\w)$ called the {\bf ball-swapping symplectomorphisms} \cite{Wu14}, which essentially captures the connecting map in the associated homotopy exact sequence.  See section \ref{BS} below for details.  Locally, Seidel's generators of the monodromies in the configuration space is exactly equivalent to swapping two balls.

As it turns out, this explicit construction has interesting applications to studying mapping classes and classifying Lagrangian submanifolds in $A_n$-surface singularities \cite{Wu14}, see also Theorem \ref{t:An}.  Later, in the authors' joint works \cite{LLW15, LLW16} with T.-J. Li, we use this construction to study $\pi_0(Symp(X,\w))$ for rational surfaces with $\chi(X)=8$ (Theorem \ref{1}). 
We expect the ball swappings construction to play a role in more questions involving a larger class of symplectic manifolds.  As a final remark, we note that the topology of the space of symplectic packings has also been studied by many researchers, see for example \cite{LP04, AL07, AP13} (dimension 4) and \cite{Hind16} (dimensions higher) just for a very incomplete list.


A purpose of this note is to give a brief account on the series of works \cite{LW12, Wu14, LLW15, LL16, LLW16}.  Another goal is to summarize some problems in this direction, most of which are folklore for a long time but have not found their due appearance in the literature.  We hope the effort of documenting these problems will help as a reference in the future.  It is instructive to point out that we have intentionally limited our scope and do not attempt to reflect all recent developments in the study of symplectomorphisms.  Especially, we have not mentioned the line of works that are closely related to mirror symmetry \cite{Ailsa14,DKK16,Siegel16,SS17,Kartal18}, etc.

  Our exposition is organized as follows.  In Section 2, we will detail the construction of ball swappings, point out the connection with Dehn twists along Lagrangian spheres, and briefly discuss the case beyond dimension 4. In Section 3, we will focus on the applications of ball-swappings to rational surfaces, and the closely related Lagrangian isotopy problems will be discussed in Section 4.  Section 5 is a immature attempt to propose a general scheme for connecting the study of space of ball-packings and the topology of configuration space in a topologist's eyes, which hopefully would inspire the seasoned experts to further the connection of this circle of ideas to higher homotopy groups of $Symp(M)$.\\

\noindent\textbf{Acknowledgements.} We are grateful to Silvia Anjos, Paul Biran, Olgta Buse, Richard Hind, Tian-Jun Li, Cheuk-Yu Mak, Dusa McDuff, Martin Pinsonnault, Leonid Polterovich, for their interests in reading a draft of this note, who provided invaluable suggestions to the authors.

\section{Ball-swapping symplectomorphisms}\label{BS}




In general, there is a way of constructing symplectomorphisms from the monodromies of moduli spaces of K\"ahler manifolds.  The ``ball-swapping" construction, which will be the main topic of our article, is a special instance of this well-known approach. In what follows, we give a brief account.  Let $E$, $B$
be smooth quasi-projective varieties, and $\pi: E \rightarrow
B$ a proper smooth morphism with a
line bundle $L \rightarrow E$ which is relatively very ample.  This means the sections of $L|E_b$ define an embedding of
$E$ into a projective bundle over $B$. Assume $M=E|_{b_0}$ is the fiber over a generic point $b_0$.  Then $M$ obtains a family of symplectic forms $\w_b$ such that the cohomology class of $\w_b$ is locally constant.  Moreover, there is a symplectic connection on the bundle $E$, taken as the symplectic orthogonal of the vertical tangent space.  Therefore, $E$ becomes a symplectic fibration, with a classifying map
  \begin{equation} \label{eq:bb}
     B \longrightarrow BSymp(M,\w_{b_0}).
  \end{equation}
 Note that this family is not always universal, but it could still yield interesting monodromies.  In particular, one may take $B$ as (part of) the coarse moduli of a certain variety $M$.  We are interested in the symplectic geometry of a smooth $M$ in this note, but it is completely justified to consider monodromies of a singular symplectic space, for example, a symplectic orbifold.

The above philosophy applied to $B=Conf^{un}_k(X)$, the configuration space of $k$ distinct unordered points on a projective variety $X$, yields a family of $M=X\#k\ov\CP^n$. Over the point $b\in B$, the fiber is the blow-up of the corresponding $k$-tuples.  One may partially compactify $B$ to a normal quasi-projective variety $\bar B$ by adding a complex codimension $1$ discriminant locus $\Delta$ where two points collide.  $E$ can be also compactified into $\bar E$, by adding the corresponding blow-ups of non-reduced points.  This gives fibers with normal singularities over $\Delta$.  Take a normal disk $D_{\bar b}$ centered at $\bar b\in \Delta$, $\bar E|_{D_{\bar b}}$ gives a Lefschetz fibration with a single Morse-Bott singularity component $\CP^{n-2}$, hence the monodromy over $\partial D_{\bar b}$ is a fibered Dehn twist along a family of isotropic $S^2$'s over this $\CP^{n-2}$ \cite{WWfamily} (see also \cite{MW15}).  This is well-known when $n=2$ and Ivan Smith pointed out to us the general case.

These simple monodromies motivate the following symplectic generalization in \cite{Wu14},
which we call the \textit{ball-swapping} (see also  \cite{LW12, BLW12}).

Suppose $X$ is a symplectic manifold.  Given two symplectic ball embeddings:

$$\iota_{0,1}:\coprod_{i=1}^n B(r_i)\rightarrow X,$$

where $\iota_0$ is isotopic to $\iota_1$ through a Hamiltonian path $\{\iota_t\}$.  From the interpretation of blow-ups
in the symplectic category \cite{MP94}, the blow-ups can be represented as

$$X^{\#\iota_j}=(X\backslash\iota_j(\coprod_{i=1}^n B_i))/\sim,\text{ for }j=0,1.$$

Here the equivalence relation $\sim$ collapses the natural $S^1$-action on $\partial B_i=S^3$.  Now assume that
$K=\iota_0(\coprod B_i)=\iota_1(\coprod B_i)$ as sets, then $\iota_t$ defines a symplectic automorphism $\widetilde{\tau}_\iota$ of
$X\backslash K$, which descends to an automorphism $\tau_\iota$ of $X^{\#\iota}:=X^{\#\iota_0}=X^{\#\iota_1}$.

\begin{defn}\label{d:BS}
   We call
$\tau_\iota$ a \textbf{ball-swapping}  on $X^{\#\iota}$.
\end{defn}

One should note that a local ball-swapping is not interesting beyond the homological level: if $X=\C^n$ and the ball-swapping has $\iota_0(B_i)=\iota_1(B_i)$ for all $i$, then one may choose an isotopy from $\{\iota_t\}$ to the identity loop, which gives an isotopy of $\tau_\iota$ to identity for $n\ge2$.

The difference between a ball-swapping symplectomorphism and the monodromies constructed from algebraic geometry is a bit subtle.  What the monodromies from the $Conf^{un}_k(X)$ capture, roughly speaking, is the topology of the configuration space given by the center of the corresponding ball-packing; while the ball-swapping a priori reflects the topology of symplectic ball-packing.  \cite{LP04} shows that these two spaces can be very different when the size of ball-embedding is large even when there is only one ball involved.  However, if we are only concerned about the symplectic mapping class group, there is evidence that the two might not be that different.  We will see more examples in the next section.

There is also a non-compact version of the ball-swapping symplectomorphism given in \cite{Wu14}.  In dimension four, the construction goes as follows.  Let $D^2\subset X^4$, and $\iota_t: B_i\hookrightarrow X$, $0\le t\le1$ be a loop of ball packing in $X$ induced by Hamiltonian in $X$ which preserves $D$.  Assume $\iota_0=\iota_1=\iota$, and that $\iota_t(B_i)\cap D$ persist to be the large disk in the respective $B_i$.  Then the ball-swapping $\tau_\iota$ constructed earlier induces a Hamiltonian $\phi^D$ in $D$.  Note that $\phi^D$ can be extended to a Hamiltonian diffeomorphism $\hat\phi^D$ in $X$ which is supported near $D$.  Then the symplectomorphism $(\hat\phi^D)^{-1}\circ\tau_\iota$ defines a compactly supported symplectomorphism in $X^\iota-\hat D$, where $\hat D$ is the proper transform after the blow-up of the ball-embedding $\iota$.  We call $(\hat\phi^D)^{-1}\circ\tau_\iota$ a \textbf{non-compact ball-swapping}.

  In essence, such a non-compact ball-swapping is given by isotoping two attaching Legendrian knots on $\partial (X-D)$, which gives a symplectomorphism that switches the corresponding handle attachments of $X^\iota-\hat D$.

  \begin{example}
    In the local model when $X=\C^2$ and $D=\C\times \{0\}$ with two ball-embedding $B_0$ and $B_1$, the symplectic completion of $X^\iota-\hat D$ is symplectomorphic to $T^*S^2$.  If one takes $\iota$ so that the restriction of $\iota_t$ to the center of $B_i$ gives a generator of $\pi_1Conf^{un}_2(D^2)$, the non-compact ball-swapping gives the Dehn twist of the zero section \cite{Wu14}.
  \end{example}

Yet another generalization of the ball-swapping construction is to replace the balls by an arbitrary open region whose boundary admits a circle action.  This again induces a symplectomorphism in the symplectic cut.

\begin{example}\label{ex:Blair}

In  \cite{RubinB}, Rubin-Blaire considers $X=\CP^3$ and the neighborhood of a genus $4$ curve, see Proposition 1.48 therein.  The idea is similar to the universal bundle construction we described in the beginning:  One constructs a nodal family of genus $4$ Riemann surfaces, which is a pair $(p: \cC \to D^2)$ where
  $\cC$ is a complex analytic subvariety of dimension $2$, and $(p: \cC \to D^2)$ is a proper holomorphic map with a single Morse singularity over $0\in D^2$.  The monodromy of this family is a Dehn twist along a separating curve. There is a sequence of inclusions
  $$ \bar{\cC} \hookrightarrow S^2\times S^2 \times D  \hookrightarrow \C P^3 \times D,$$ which respect the bundle structure over the disc.  The blow-up of $\cC$ in $\C P^3 \times D$ yields a new monodromy in a Fano $3$-fold $\hat X$, which is ``ball-swapping'' type symplectomorphism from a path of embedded symplectic curves.  It should be pointed out that Rubin-Blair's construction and study of this example involves many non-trivial applications of algebro-geometric considerations, which enables one to extract information of the loop of the genus $4$ curve.


\end{example}

\section{Symplectic mapping class groups of rational surfaces} 
\label{sub:symplectic_mapping_class_groups_of_rational_surfaces}



One of the most interesting applications for the ball-swapping constructions is to study the symplectomorphism groups of blow-ups when $X=\CP^2$.  From a theorem due to McDuff \cite{McD96}, the space of ball-packing is connected.  Denote
\begin{equation}
\begin{aligned}
  &Emb(c_i):=\\
  &\{(f(\coprod_{i=1}^k B_i(c_i)),p_i)| f:\coprod_{i=1}^k B_i(c_i)\hookrightarrow \CP^2
  \text{ is an embedding
  of symplectic open balls, }\\
  &p_i\in B_i\}\end{aligned}
\end{equation}

as the space of (unparameterized) ball-packings in $\CP^2$.  Note that the extra marked points $p_i$ will allow the evaluation map to be defined, and the forgetful map that removes them induces a homotopy equivalence to the ordinary ball-packings.  Throughout this section, we only consider the case when $Emb(c_i)$ is non-empty, which boils down to a concrete numerical condition in $c_i$ involving Cremona transforms.  See \cite{MS12}, or the reduction algorithm in \cite[Section 2.3]{BP13}. From a homotopy theoretic point of view, the ball-swapping is a consequence of the following action-orbit fibration

\begin{equation}\label{e:action}
     Ham(\CP^2;\coprod_{i=1}^k B_i)\to Ham(\CP^2)\to Emb(c_i).
\end{equation}

Here, $Ham(\CP^2;\coprod_{i=1}^k B_i)$ denotes the group of Hamiltonian symplectomorphisms on $\CP^2$ that preserves the image of the embeddings of $\coprod_{i=1}^k B_i$ (but we do not require the image of individual $B_i$'s being preserved).  We note that even though $Ham(\CP^2)$ is connected, $Ham(\CP^2;\coprod_{i=1}^k B_i)$ is usually disconnected.  On the other hand, there is a natural map
$$\alpha: Ham(\CP^2;\coprod_i B_i)\to Symp(\CP^2\#k\ov\CP^n).$$
$\alpha$ is a homotopy equivalence if one assumes that the exceptional curves obtained by blowing up $B_i$ have the minimal symplectic area (see \cite[Theorem 2.5]{LP04}).   This holds, for example, when

\begin{equation}\label{e:condition}
     \w(E_i)=\frac{1}{m}\w(H)\text{ for }m\in \Z, \forall i.
\end{equation}

We also note that, if we denote $Ham^0(\CP^2;\coprod_{i=1}^k B_i)<Ham(\CP^2;\coprod_{i=1}^k B_i)$ to be the subgroup that preserves individual images of $B_i$, then 

$$\alpha_0:Ham^0(\CP^2;\coprod_{i=1}^k B_i)\to Symp_h(\CP^2\#k\ov\CP^n)$$

also yields a homotopy equivalence with the same non-degenerating assumption (recall from the introduction $Symp_h(\CP^2\#k\ov\CP^n)$ is symplectomorphism group that acts trivially on homology). 
Without \eqref{e:condition}, $\alpha$ and $\alpha_0$ at least induce isomorphisms on $\pi_0$, because the bubbling phenomenon is of codimension $2$ (see the argument in \cite{LP04}).  The induced homotopy sequence therefore yields (for any packing sizes $c_i$),

$$\pi_1(Ham(\CP^2))\cong \Z_3\to \pi_1(Emb(c_i))\to \pi_0(Symp(\CP^2\#k\ov\CP^n))\to 1.$$

It should be clear from our construction that the ball-swappings are those symplectomorphisms in the image of $\pi_1(Emb(c_i))$, with the kernel being an order 3 ``full twist".  This implies

\begin{lemma}\label{l:AllBS}
    Any symplectomorphisms in $Symp(\CP^2\#k\ov\CP^2)$ is Hamiltonian isotopic to a ball-swapping.
\end{lemma}

Given the similarity of ball-swappings and the monodromies from algebraic geometry, it is natural to make the following conjecture, which has been open for decades, and is at least implicit from Seidel's earlier works \cite{Se4dim}.

\begin{conj}\label{conj:strong}
  $\pi_0Symp(\CP^2\#k\ov\CP^2)$ is generated by Lagrangian Dehn twists.
\end{conj}

When $k\le 4$, this was confirmed in \cite{LLW15}, following the work of \cite{Ev11}.  When $k=5$, $\pi_0Symp(\CP^2\#k\ov\CP^2)$ were computed in \cite{Se4dim}, \cite{Ev11} and \cite{LLW16}.

\begin{thm}\label{1}
Let $  X $ be  $ \C P^2  \# n{\overline {\C P^2}}.$
If $ n\leq 4,$   $\pi_0 (Symp_h (X,\w))$ is trivial, and hence  $\pi_0 (Symp (X,\w))$ is a finite reflection group. When $n=5$, the moduli space of symplectic form up to diffeomrophism and rescalling has real dimension 5, and  $\pi_0 (Symp_h (X,\w))$ is completely known to be: \\
$\bullet$ the 5-strand pure braid group of $S^2$ on the monotone point (\cite{Se4dim, Eva11}); \\
$\bullet$ the 4-strand pure braid group of $S^2$ on a one-dimensional family starting at the monotone point (\cite{LLW16});\\
$\bullet$ trivial in every other cases (\cite{LLW16}).
\end{thm}

  Although the above works did not identify the generators as Lagrangian Dehn twists, the expected set of Lagrangian spheres can be found by \cite{LW12}, and the conjecture seems likely true.

While all these known instances are in favor of Conjecture \ref{conj:strong}, it remains largely unsolved at the time of writing.  From the classification of Lagrangian sphere classes \cite{LW12}, the existence of Lagrangian spheres depends on countable numerical constraints on the class $[\w]$.  Therefore, a generic choice of symplectic form does not admit a Lagrangian sphere.  Therefore, Conjecture \ref{conj:strong} implies the following weaker version, which is still open.
\begin{conj}\label{conj:weak}
  For generic choice of symplectic form $\w$, $Symp(\CP^2\#k\ov\CP^2)$ is connected.
\end{conj}

As another piece of evidence for the above conjectures, \cite{LLW16} showed that when $k=5$, Conjecture \ref{conj:weak} holds.

When $k\le 8$, the symplectic rational surfaces admit monotone symplectic forms, yielding a series of \textit{symplectic Del Pezzo surfaces} $X_k$, $k\le 8$.  When $k\le5$, nice descriptions of the homotopy type of $Symp(X_k)$ has been obtained, summarized as follows

\begin{itemize}
  \item $Symp(X_0)\sim PU(3)$, $Symp(X_1)\sim U(2)$ (Gromov \cite{Gro85}).
  \item $Symp(X_2)\sim T^2\rtimes \Z/2$ (Lalonde-Pinsonnault \cite{LP04}).
  \item $Symp_h(X_3)\sim T^2$, $Symp_h(X_4)\sim pt$, $Symp(X_5)\sim \Diff^+(S^2,5)$ (Evans, \cite{Ev11}).
\end{itemize}

Here, $\sim$ denotes weakly homotopy equivalences. It is natural to wonder

\begin{prob}\label{p:closedSMC}
  What are the homotopy type for $Symp(X_k)$ when $k=6,7,8$?
\end{prob}

This problem was also first considered by Seidel \cite{Se4dim} where he briefly reasons why one should expect the answer to be related to a braid group on a higher genus curve from the moduli space point of view.

From a more technical standpoint, the proof of all known cases depends on the choice of a nice curve configuration in the corresponding Del Pezzo surface.  While these configurations for $k\le5$ are all simply connected, $X_6$, $X_7$ and $X_8$ seems to require one to consider certain non-simply connected curve configurations, which again leads us to certain braid groups on higher genus curves.  See \cite[Example 1.13]{Se4dim} for a more detailed discussion.\\

Conjecture \ref{conj:strong} and Problem \ref{p:closedSMC} also have an open counterpart.  From the study of Fukaya categories, Khovanov and Seidel \cite{KS02} made the following deep observation.

\begin{thm}[\cite{KS02}]\label{t:KS}
  There is a braid group embedding $Br_{n+1}\hookrightarrow \pi_0Symp_c(W_n^m)$, where $W_n^m$ is the compactly supported symplectic mapping class group of the $A_n$-Milnor fiber of complex dimension $m$.
\end{thm}

Here, $W_n^m$ is the plumbing of $n$ copies of $T^*S^m$, and the braid group generators are given by Lagrangian Dehn twists of zero sections.  In general, we know very little about the higher dimensional symplectomorphism groups (even the connectedness of $Symp_c(B^{2n})$ for $n>2$ is widely open), but the following problem seems in reach of the current technology.


\begin{prob}\label{p:ADE}
   Are the symplectic mapping class groups of $ADE$-type Milnor fibers generated by Lagrangian Dehn twists?  And what are $\pi_0Symp_c(W_n^2)$?
\end{prob}

So far, we have a complete answer to this question for $A_n$-Milnor fibers in dimension $4$.

\begin{thm}[\cite{Ev11,Wu14}]\label{t:An}
   The compactly supported symplectomorphism group $Symp_c(W^2_n)$ is weakly homotopy equivalent to the $(n+1)$-strand braid group on the disk $Br_{n+1}$.  Moreover, the Dehn twists of the plumbed zero sections corresponds to the generators of the braid group.
\end{thm}

Evans first showed in \cite{Ev11} that $Symp_c(W_n)$ surjects onto $Br_{n+1}$.  His proof again reduces the problem to a compactified rational surface and computes the long exact sequence of homotopy groups by the action-orbit fibration for certain configuration of symplectic curves.  The second author constructed the generators of the braid group using the non-compact ball-swappings and proved that they are Hamiltonian isotopic to the Lagrangian Dehn twists along the zero sections.

The main difficulty of Problem \ref{p:ADE} is as follows. To turn all symplectomorphisms of $A_n$-Milnor fibers into a non-compact ball-swapping, a very concrete model for the compactification is used in \cite{Wu14}.  However, compactifications of $D$ and $E$ type Minor fibers do not have a similar interpretation, and certain generalization of the ball-swapping symplectomorphisms might play a role in this problem.  A much weaker statement was obtained by Mak and the second author \cite{MW15}, which at least provides a piece of evidence for the positive side of Problem \ref{p:ADE}.

\begin{thm}\label{t:SympWgen}
    Let $W$ be an $ADE$-Milnor fiber of any dimension.  For any compactly supported symplectomorphism $\phi\in Symp_c(W)$, the corresponding auto-equivalences $\Phi_\phi\in D^\pi Aut(\fuk(W))$  is split generated by the corresponding auto-equivalences induced by compositions of Dehn twists along the standard vanishing cycles.
\end{thm}

For more general tree-like plumbings of $T^*S^2$ the same problem is still interesting but even less transparent.  In another direction, one may consider non-contractible loops in $D$ in the non-compact ball-swapping construction.

\begin{prob}
  Suppose $D^2\subset X^4$ is dual to a multiple of $[\w]$.  Let $\iota_t(B(c))$ be a family of ball-embeddings satisfying the assumption in the construction of non-compact ball-swappings, and assume $\iota_t(p)$, the trace of the center of $B(c)$, is a non-contractible loop in $D$.  Is the resulting non-compact ball-swapping always non-trivial in the compactly supported mapping class group?
\end{prob}

\begin{example}
  Let $D:=\cup_{i=1}^3l_i\subset \CP^2$ be the union of three generic lines, and choose a path of ball embeddings whose center traces the non-contractible loop in $l_1-(l_2\cup l_3)$, which is a $\C^*$.  Although in this case, $D$ is not smooth, the construction still yields a compactly supported symplectomorphism in $T^*T^2$ with a two handle attached.

  More concretely, take the $\partial^\infty T^*T^2$ as the standard $T^3$, the handle is attached along a Legendrian representing the meridian of $l_1$.  Roughly, the ball-swapping represents an isotopy of this Legendrian, which projects onto a non-contractible circle of the Lagrangian $T^2$.  Unfortunately, we do not have a good description in this handle attachment picture for the ``rollback'' $\hat\phi^D$ in the construction of the non-compact ball-swapping, which is the key to getting a compactly supported symplectomorphism.

  When $D=Q\cup l\subset CP^2$, the union of a generic smooth quadric and a line, its complement in $\CP^2$ is the line complement of $T^*\RP^2$, which is a well-known example in mirror symmetry.  Take a loop on the quadric which divides its area in half, then there is also a corresponding ball-swapping along this path.

  It would be interesting to know whether these symplectomorphisms are non-trivial in the compactly supported symplectomorphism groups.
\end{example}

We end the discussion on $\pi_0(Symp(X))$ by remarking that, there should be a higher dimensional analog of the non-compact ball-swapping, which replaces the movement of ball-packing by an isotopy of an Legendrian knot.  The harder part is to find an appropriate setting to allow the ``rollback map'' at infinity to obtain a compactly supported symplectomorphism.\\





In the rest of this section, we discuss briefly the fundamental group of $Ham(X)$ for a rational surface.
In \cite{McD08}, Dusa McDuff gives an approach to obtain the upper bound of
$\pi_1Ham(X,\omega)$, where $(X, \omega)$ is a symplectic rational 4 manifold.
 We can follow the route of Proposition 6.4 in \cite{McD08} to give
a proof of the following result, refining \cite[Proposition 6.4, Corollary 6.9]{McD08}.  See also \cite[Proposition 4.13]{LL16}.

\begin{thm}\label{l:mcduff64}
   Let $(X,\w)$ be a symplectic rational surface with $b_2(X)=r$, and $(\wt X,\wt\w)$ be the blow-up of $X$ for $k$ times.  If all the new exceptional divisor $E_i$ has equal area, and this area is strictly smaller than all exceptional divisors in $X$,
   $$rank[\pi_1(Symp_h(\widetilde{X}))]\le rank[\pi_1(Symp_h({X}))]+kr.$$

    If the new exceptional divisor  smaller than all exceptional divisors in $X$, then
   $$rank[\pi_1(Symp_h(\widetilde{X}))]\le \pi_1(Symp_h({X}))+r.$$

\end{thm}

In addition, there is a lower bound of the rank of $\pi_1(Ham(X,\w))$ given by counting circle actions in \cite{Ped15}. Note that this lower bound may not be optimal. In \cite{LL16} \cite{LLW16} and \cite{LLW3}, a finer lower bound  is given by considering the long exact sequence of the fibration with total space being $Symp_h(X,\w)$ and taking into account of $\pi_0(Symp_h(X))$.

Note that in most cases, this lower bound agrees with the upper bound given by Theorem \ref{l:mcduff64} and hence is the precise rank of $\pi_1(Symp_h(X,\w))$.  Define $Q(X):=\frac{1}{2}(\chi(X)-2)(\chi(X)-3)$.  Then the above result can be summarized for rational surfaces $(X,\w)$ with $\chi(X)\le7$, where $\pi_0(Symp_h(X,\w))$ is trivial, as $$Q(X)=r^+[\pi_0(Symp(X,\omega))]+  Rank [ \pi_1(Symp_h(X,\omega))].$$
here $r^+[\pi_0(Symp(X,\omega))]$ is the number of homology classes admitting $(-2)$ Lagrangian spheres.  These classes only depends on the cohomology class of $\w$ by a theorem in \cite{LW12}.

When $8\le \chi(X)\le 11$, due to the fact that $\pi_0(Symp_h(X,\w))$ could be non-trivial, we consider the following:

\begin{conj}[Li-Li, \cite{LL16}]
    The rank of $\pi_1(Symp_h(X,\w))$ and $\pi_0(Symp_h(W,\w))$ satisfies  $$Q(X)=r^+([\pi_0(Symp(X,\omega))])+  Rank [ \pi_1(Symp_h(X,\omega))]-Rank [\pi_0(Symp_h(X,\omega))].$$

\end{conj}

This conjecture was verified in  \cite{LLW16} when $\chi(X)=8$, where we computed that that  $Q(X, \w)=15$ when $\chi(X)=8$.  Also note that this conjecture provides information for Problem \ref{p:closedSMC} above.

Note that the result of Anjos-Pinsonnault \cite{AP13} and Anjos-Eden \cite{AE17} also imply the above computation of $\pi_1(Symp_h(X,\w))$ for any given form on 3-fold blow up of $\C P^2$.  They give a generating set of $\pi_1(Symp_h(X,\w))$ using
circle action. In contrast, our approach in \cite{LLW16} gives a generating set of $\pi_1(Symp_h(X,\w))$ using $(-2)$ symplectic spheres.  The two generating sets have the following correspondence.
For any $(-2)$ symplectic sphere in 3 fold blow-up of $\C P^2$, there
  is a semi-free circle $\tau$ action fixing it, where $\tau$ is a generator of $\pi_1(Symp_h(X,\w))$.  One may then bring up the following natural question: can generators of $\pi_1(Symp_h(X,\w))$ always be represented by Hamiltonian or symplectic circle actions for toric surfaces? This is closely related by Seidel representation \cite{Sei99} with Quantum cohomology, see \cite{MT06} and \cite{Gon06} for related works.

\section{Lagrangian isotopy classes in dimension four} 
\label{sec:lagrangian_homotopy_classes_in_dimension_}

  The question of symplectic mapping class groups is closely tied to the Lagrangian isotopy problems.  Recall that a \emph{Biran decomposition} of a symplectic $4$-manifold $X$ with $[\w]\in H^2(X,\Z)$ is a decomposition $X=\Lambda\coprod B$, where $B$ is a disk bundle over a symplectic Donaldson hypersurface $\Sigma$, and $\Lambda$ is an \textit{isotropic skeleton} given by the limit of an inverse Liouville flow of $X-\Sigma$ \cite{Bi01}.

In \cite{Coffey05}, Coffey showed the following result.

\begin{thm}\label{t:Coffey}
       Let $(X, \w)$ be a symplectic $4$-manifold with Biran decomposition $\Lambda\coprod B$, such that $\Lambda$ is a smooth Lagrangian submanifold and $\Sigma$ is of genus zero. Then $Symp(X)$ is homotopy equivalent to the $\eL$, the space of Lagrangian embedding of $\Lambda$.
 \end{thm}

 There is also a version for non-smooth $\Lambda$, for which the statement is considerably more involved, for which we refer the reader to \cite[Theorem 2.6]{Coffey05}.

One should note that the general problem of classifying Lagrangian embeddings is extremely difficult. In \cite{Se99}, Seidel constructed homologous but non-hamiltonian Lagrangian spheres; and in \cite{Vid06}, Vidussi showed that certain homology classes admit a class of Lagrangian tori that are mutually smoothly non-isotopic.  

We focus on the classification up to Hamiltonian isotopy.  By Weinstein's theorem, this problem is more difficult than the famous Arnold's nearby Lagrangian conjecture, which asserts all exact Lagrangians in a cotangent bundle are all Hamiltonian isotopic to the zero section.  So far, this strongest version of nearby Lagrangian conjecture is only known in the case of $T^*S^2$, $T^*\RP^2$ and $T^*T^2$, see \cite{Hind04} and \cite{DRGI16}.  Through the case of $T^*S^2$ and $T^*\RP^2$, we also know that the space of exact Lagrangian embeddings is contractible by \cite{HPW13}.  The nearby Lagrangian conjecture for the Klein bottle and the classification of Lagrangian isotopy classes in $\C^2$ seems to lie within the reach of current technology, but the situation for $T^*\Sigma_g$ for $g>1$ might require new insights.  In this direction, Hind and Ivrii \cite{HindIvrii2} showed the smooth isotopiness of all Lagrangian embeddings.

For the reasons stated above, we will mainly focus on the Lagrangian isotopy problems of $S^2$ and $\RP^2$.  Hind first investigated the Hamiltonian uniqueness of Lagrangian spheres in $S^2\times S^2$ \cite{Hind04}. In the joint works of Borman, Li and the second author, it is shown that for all Lagrangian $S^2$'s are smoothly isotopic in a symplectic rational or ruled surface.  On the other hand, they are unique up to symplectomorphisms \cite{LW12, BLW12}.    Moreover, it was proved that homologous Lagrangian spheres are Hamiltonian isotopic to each other in $\CP^2\#k\ov\CP^2$ when $k\le 4$, and Lagrangian $\RP^2$ are Hamiltonian isotopic to the standard real locus in $\CP^2$ in \cite{LW12}.  

One should note that the smooth uniqueness of homologous Lagrangian spheres also follows from Shevchishin's work \cite[Corollary 2(ii)]{She10}, combined with the uniqueness up to symplectomorphisms proved in \cite{BLW12}.  Recently, Shevchishin and Smirnov announced a new result on classifying Hamiltonian isotopy classes of Lagrangian spheres in $\CP^2\#5\ov\CP^2$.  It is certainly an interesting question to ask for a similar classification for $k>5$, but it likely depends on how much we can extract from Problem \ref{p:closedSMC}.

Interestingly, the following question remains unsolved at the moment.

\begin{prob}\label{p:RP2}
   Are $\Z/2$-homologous Lagrangian $\RP^2$ unique up to symplectomorphisms and smooth isotopies in symplectic rational surfaces?
\end{prob}

To explain the context, we recall that in \cite{BLW12}, the treatments of Lagrangian $\RP^2$ and $S^2$ are largely similar.  First, we perform a symplectic cut near a given Lagrangian $L=\RP^2$ or $S^2$, which results in a symplectic $(-4)$ or $(-2)$-sphere.  We call a set of exceptional curves \textit{push off curves} if they are mutually disjoint and disjoint from the given $(-4)$ or $(-2)$-sphere, and that the blow down of them gives a symplectic $S^2\times S^2$.  A neck-stretching argument will show the existence of push off curves for any $(-4)$ or $(-2)$ sphere.  Now take a homologous Lagrangian $L'$ and proceed with the same process, one can find another complete set of exceptional spheres.

The story diverges from here.  In the case of $(-2)$-spheres (corresponding to Lagrangian $S^2$'s), the two sets of disjoint exceptional curves are homologous, hence are Hamiltonian isotopic to each other.  At this point, we use a symplectic sum surgery to replace the symplectic $(-2)$ sphere by the original Lagrangian $L$ and $L'$.  This finds, up to a Hamiltonian displacement of $L'$, a set of exceptional spheres that are disjoint simultaneously from $L$ and $L'$, whose blow-down will yield a symplectic $S^2\times S^2$.  Then the smooth isotopy of $L$ and $L'$ will follow from blowing down this set of exceptional curves and reducing the problem to the well-known case of $S^2\times S^2$.  To get the uniqueness up to symplectomorphisms, blow-down the same push-off curves.  If one traces McDuff's proof of connectedness of space of ball-packing, it is not hard to be convinced a slight modification will show that the same argument works in the complement of a $(-2)$-sphere.  Therefore, the isotopy of these balls will result in a ball-swapping, which again yields a symplectomorphism that sends $L$ to $L'$.

The main issue of Problem \ref{p:RP2} lies in that, $\RP^2$ does not have a $\Z$-homology class but only a $\Z/2$-class.  Therefore, the two complete sets of exceptional spheres disjoint from two distinct $\Z/2$-homologous Lagrangian $\RP^2$'s can have different homology classes.  This is equivalent to the existence of two smooth complete sets of exceptional curves that are $\Z/2$-homologous but not $\Z$-homologous.  A concrete example was given in \cite{BLW12} where this actually happens, which results in two smoothly non-isotopic but $\Z/2$-homologous Lagrangian $\RP^2$'s when $k=10$ if we allow the symplectic form to vary.

However, it is noted in \cite{BLW12} that this is the only possible obstruction for Problem \ref{p:RP2}.  In other words, if the complete sets of exceptional spheres are homologous, the above proof for Lagrangian $S^2$ can be carried over to the Lagrangian $\RP^2$ without any difficulties.   This is the case when $k\le8$, where Borman-Li-Wu concluded the smooth isotopy uniqueness for Lagrangian $\RP^2$'s.  On the other hand, given two complete sets of exceptional classes that are $\Z/2$-homologous but not $\Z$-homologous, we have not found a symplectic form on any rational surfaces such that both of these sets can be made symplectic and disjoint from a fixed symplectic $(-4)$-sphere, due to a symplectic packing obstructions.

\begin{example}
   Take $\CP^2\#10\ov\CP^2$.  Lemma 4.12 of \cite{BLW12} constructs $L$, a Lagrangian $\RP^2$ for certain symplectic form $\w$, which has the $\Z/2$-class $E_1+E_2+E_3$.  It is shown by Lemma 4.11 in \cite{BLW12} that, this $\RP^2$ can not be smoothly isotoped away from the exceptional curve in class $3H-2E_1-E_4-E_5-\cdots-E_9$, but their pairing clearly vanishes in $\Z/2$-homology.

   On the other hand, there is a smooth diffeomorphism $\phi$ which exchanges the homology class of $E_{10}$ and $3H-2E_1-E_4-E_5-\cdots-E_9$.  This can be constructed as a smooth Dehn twist along a smooth $(-2)$-sphere in class $3H-2E_1-E_4-E_5-\cdots-E_{10}$.  Therefore, $\phi(L)$ is a $\phi_*\w_0$-Lagrangian $\RP^2$, which is not smoothly isotopic but $\Z/2$-homologous to $L$.

   However, we cannot find a way to upgrade $\phi$ to a symplectomorphism.  This upgrade is equivalent to find a ball packing $B_i(c_i)$ in $\CP^2$, $i=1,\cdots,10$, such that $E_1+E_2+E_3$ admits a Lagrangian $\RP^2$ embedding, and $\w(3H-2E_1-E_4-E_5-\cdots-E_{10})=0$.  These conditions can be translated into very concrete numerical conditions from a Cremona criterions derived from \cite{LL02}, see \cite{MS12} \cite{Pin08emb} or \cite{BP13}, and a very recent development of packing of $\RP^2$ \cite{ShSm19}.

   If one puts together these numerical constraints, the problem becomes an elementary Lagrangian multiplier problem.  The upshot is, such a packing doesn't exist (even without taking Shevchishin-Smirnov's triangle inequality \cite{ShSm19} into account).  Unfortunately, the authors are not aware of an effective algorithm to exploit these numerical constraints to solve Problem \ref{p:RP2} in complete generality. 
\end{example}

Given what we knew about Lagrangian spheres in rational surfaces, as well as the symplectic mapping class groups in Milnor fibers, it is natural to expect the following.

\begin{conj}
  Lagrangian spheres in surface Milnor fibers of $ADE$ type singularities are Hamiltonian isotopic to one of the zero sections after a sequence of Dehn twists along these standard spheres.
\end{conj}

The conjecture is known to be true for the $A_n$ type surface Milnor fibers from the characterization of their mapping class groups, see \cite{Wu14}.  It would also be interesting to study the corresponding questions in the cusp singularities \cite{Kea15}

We conclude this section by raising a general question regarding the displacement energy.  Chekanov discovered in \cite{Chek98} that, a displaceable Lagrangian $L$ can be displaced by a Hamiltonian diffeomorphism $\phi$ only if the Hofer energy of $\phi$ is larger than the smallest holomorphic disk with boundary on $L$.  There are other interesting conjugation invariant norms in $Ham(X)$.  Recall that, one might denote as the \textit{autonomous norm} as

$$||f||_{aut}:=min\{k|f_1\cdot f_2\cdots f_k=f, \text{ where $f_i$ are autonomous} \},$$

see \cite{BK13} for example.  There are a series of interesting recent developments in the study of autonomous norms \cite{PS16, ZhangJun}, which proves autonomous Hamiltonians are only a small subset among all Hamiltonian diffeomorphisms in general.  Therefore, it is natural to wonder:

\begin{prob}
   Is there an example of Lagrangian submanifold $L$ such that $Disp_{aut}(L)=k$ for every $k\in \Z^+$?
\end{prob}

Here, $Disp_{aut}(L)$ denotes the smallest autonomous norm for a Hamiltonian to displace $L$.  Unfortunately, as of the time of writing, it is even unclear if one could find a Lagrangian so that $k>0$, that is, a Lagrangian which can be displaced, but not displaceable by an autonomous Hamiltonian.  It is not clear if the ball-swappings that are isotopic to identity potentially gives such examples.

\section{The topology of configuration spaces: from classical to symplectic} 
\label{sec:floer_theory_of_ball_swappings}




As is already clear from our previous discussions, ball-swappings are closely related to the configuration spaces of points of the ambient symplectic manifold.  However, if one focuses on the symplectomorphism group of blow-ups, an additional layer of complication adds up.  For instance, when we consider the unordered configurations of more than 3 points in $\CP^2$, one has to also include the discriminant locus where three points are collinear.  While this is very interesting in the algebraic geometric point of view, the following problem makes a more direct connection to a topologist's configuration space.

\begin{prob}\label{p:classical}
  Given any symplectic manifold $X$, consider the forgetful map $$\sigma:Emb(c_i)\to \Conf_k(X)$$ where $\sigma$ restricts a ball-packing to the map of the marked point $p_i$.  What is the  kernel/cokernel of the induced homomorphism on cohomology $H^i(\sigma)$?
\end{prob}

The topology of configuration spaces is an active area of research with a long history, just to name very few examples, see \cite{McD75,
Arnold69, Church12}.
A wishful speculation in the direction of Problem \ref{p:classical} is a higher homotopical generalization of Biran's stabilization theorem

\begin{conj}\label{conj:configuration}
  Assume $X$ has dimension $4$ in Problem \ref{p:classical}.  then for any $m\in\Z^+$, there exists an $\epsilon(m)>0$, when $c_i<\frac{\epsilon(m)}{k}$, such that $\sigma$ is  $m$-connected.
\end{conj}

For general sizes of blow-ups, one cannot expect the homotopy equivalence to hold.  Instead, the space of ball-packing should exhibit more interesting symplectic behaviors, and the discrepancy between classical and symplectic geometry is partly captured by the forgetful map.

In the rest of our discussions, we will only consider cohomology with $\Q$-coefficients.
Pioneering results in the direction of Conjecture \ref{conj:configuration} were obtained by Anjos, Lalonde and Pinsonnault \cite{LP04} \cite{Pin08} \cite{ALP} \cite{AP13}, where they studied one-ball packing in $S^2$ bundle over $S^2$ and up to two-ball packing in $\CP^2$.


The main strategy of the computation of cohomology rings by Lalonde and Pinsonnault relies on a computation of the rational homotopy group of $Symp(S^2\times S^2\#\ov\CP^2)$ based on \eqref{e:action} and the map $\alpha$.  To get to the discussion of general ball-packing sizes, the key point is to analyze the holomorphic curves of negative squares and study the topology of the space they form.  With all these hard work, note that $Symp(S^2\times S^2\#\ov\CP^2)$ is a topological group, which in particular is an $H$-space.  There is then a connection between the rational homotopy and the cohomology ring by a theorem due to Cartan-Serre.  From this, one should be able to deduce the rational homotopy type of $Emb(c_i)$ and even the homotopy type in nice situations.  However, to relate the cohomology of the space of ball-packing to $Conf^{ord}_k(\CP^n)$ via this approach factors through the homotopy group of $Symp(X\#k\ov\CP^2)$, and almost always gets into the homotopy level.  This might not always be easy.  Floer theory could be an alternative tool to attack this problem, and one of the closest related references could be Varolgunes's thesis \cite{Varothesis}, from which one can obtain a symplectic cohomology for the complement of the ball-packing.

Another potentially interesting relation with classical topology is the study of loop space of ball-packings.  Assuming $X=\CP^2$ and $c_i$ satisfy \eqref{e:condition}, we have homotopy equivalence $\alpha$ in Section \ref{sub:symplectic_mapping_class_groups_of_rational_surfaces}.  Therefore, \eqref{e:action} is equivalent to
\begin{equation}\label{e:action2}
     Symp(\CP^2\# k\ov\CP^2)\to Ham(X)\to Emb(c_i).
\end{equation}

 Replace \eqref{e:action2} by an equivalent homotopy fibration
\begin{equation}\label{e:loopfibration}
     \Omega Emb(c_i)\xrightarrow{\varphi} Symp(\CP^2\# k\ov\CP^2)\to Ham(X).
\end{equation}

It is shown by Xicotencatl \cite{Xic02} that 
$$\Omega \Conf_k(\CP^2)\sim (S^1)^k\times\Omega F_{S^1}(S^5),$$ 
where $F_{S^1}(S^5)=\{(z_1,\cdots,z_k)\in (S^5)^k: z_iS^1\cap z_jS^1=\emptyset\text{ if }i\neq j\}$, where the $S^1$ acts by the standard Hopf fibration.  This yields a distinguished $(S^1)^k$ factor in the homotopy type of $\Omega\Conf_k(\CP^2)$.  As of now, it is unclear how much the classical topology of loop space of configurations contributes to $Symp(\CP^2\# k\ov\CP^2)$, but from the current computations in the literature, we make the following speculation.

\begin{conj}
   Consider the looped forgetful map $\Omega\sigma:\Omega Emb(c_i)\to \Omega Conf_k(\CP^2,k)$.  Then $\varphi^*H^1(Symp(\CP^2\#k\ov\CP^2))\cap (\Omega\sigma)^*(H^1((S^1)^k))=0$.
\end{conj}

In other words, we expect that the $(S^1)^k$-factor in the loop space of the configuration space does not contribute to the blow-up of the symplectomorphism groups.  In general, we expect the complication of $Symp(\CP^2\#k\ov\CP^2)$ comes mostly from the discrepancy between $\Omega Emb(c_i)$ and $\Omega \Conf_k(\CP^2)$.

\bibliography{FukRef}

\newcommand{\etalchar}[1]{$^{#1}$}
\begin{thebibliography}{CCGF{\etalchar{+}}14}

\bibitem[Abr98]{Abr98}
Miguel Abreu.
\newblock Topology of symplectomorphism groups of {$S^2\times S^2$}.
\newblock {\em Inventiones Mathematicae}, 131:1--23, 1998.

\bibitem[AE19]{AE17}
Silvia Anjos and Sinan Eden.
\newblock The homotopy {L}ie algebra of symplectomorphism groups of 3-fold
  blow-ups of $s^2\times s^2, \omega_{std}\bigoplus \omega_{std}$.
\newblock {\em Michigan Math Journal, Advance publication}, 2019.

\bibitem[AGK09]{AGK09}
Miguel Abreu, Gustavo Granja, and Nitu Kitchloo.
\newblock Compatible complex structures on symplectic rational ruled surfaces.
\newblock {\em Duke Math. J.}, 148:539--600, 2009.

\bibitem[AL07]{AL07}
Slvia Anjos and Fran\c{c}ois Lalonde.
\newblock The topology of the space of symplectic balls in {$S^2\times S^2$}.
\newblock {\em C. R. Math. Acad. Sci. Paris}, 345(11):639--642, 2007.

\bibitem[ALP09]{ALP}
Silvia Anjos, Fran\c{c}ois Lalonde, and Martin Pinsonnault.
\newblock The homotopy type of the space of symplectic balls in rational ruled
  4-manifolds.
\newblock {\em Geom. Topol.}, 13(2):1177--1227, 2009.

\bibitem[AM00]{AM00}
Miguel Abreu and Dusa McDuff.
\newblock Topology of symplectomorphism groups of rational ruled surfaces.
\newblock {\em J. Amer. Math. Soc.}, 13(4):971--1009 (electronic), 2000.

\bibitem[Anj02]{Anj02}
Silvia Anjos.
\newblock Homotopy type of symplectomorphism groups of {$S^2\times S^2$},.
\newblock {\em Geom. Topol}, pages 195--218, 2002.

\bibitem[AP13]{AP13}
Silvia Anjos and Martin Pinsonnault.
\newblock The homotopy {L}ie algebra of symplectomorphism groups of 3-fold
  blow-ups of the projective plane.
\newblock {\em Math. Z.}, 275(1-2):245--292, 2013.

\bibitem[Arn69]{Arnold69}
V.~I. Arnold.
\newblock The cohomology ring of the group of dyed braids.
\newblock {\em Mat. Zametki}, 5:227--231, 1969.

\bibitem[BH11]{BH11}
Olguta Buse and Richard Hind.
\newblock Symplectic embeddings of ellipsoids in dimension greater than four.
\newblock {\em Geom. Topol.}, 15(4):2091--2110, 2011.

\bibitem[BH13]{BH13}
O.~Buse and R.~Hind.
\newblock Ellipsoid embeddings and symplectic packing stability.
\newblock {\em Compos. Math.}, 149(5):889--902, 2013.

\bibitem[Bir96]{Bi96}
Paul Biran.
\newblock Connectedness of spaces of symplectic embeddings.
\newblock {\em Internat. Math. Res. Notices}, (10):487--491, 1996.

\bibitem[Bir97a]{Bi97}
P.~Biran.
\newblock Symplectic packing in dimension {$4$}.
\newblock {\em Geom. Funct. Anal.}, 7(3):420--437, 1997.

\bibitem[Bir97b]{Bir97}
Paul Biran.
\newblock Geometry of symplectic packing.
\newblock PhD Thesis, Tel Aviv University, 1997.

\bibitem[Bir01a]{Bi01}
P.~Biran.
\newblock Lagrangian barriers and symplectic embeddings.
\newblock {\em Geom. Funct. Anal.}, 11(3):407--464, 2001.

\bibitem[Bir01b]{Bir01}
Paul Biran.
\newblock From symplectic packing to algebraic geometry and back.
\newblock In {\em European {C}ongress of {M}athematics, {V}ol. {II}
  ({B}arcelona, 2000)}, volume 202 of {\em Progr. Math.}, pages 507--524.
  Birkh\"{a}user, Basel, 2001.

\bibitem[BK13]{BK13}
Michael Brandenbursky and Jarek Kedra.
\newblock On the autonomous metric on the group of area-preserving
  diffeomorphisms of the 2-disc.
\newblock {\em Algebr. Geom. Topol.}, 13(2):795--816, 2013.

\bibitem[BLW14]{BLW12}
Matthew~Strom Borman, Tian-Jun Li, and Weiwei Wu.
\newblock Spherical lagrangians via ball packings and symplectic cutting.
\newblock {\em Selecta Mathematica}, 20(1):261--283, 2014.

\bibitem[BP13]{BP13}
Olguta Buse and Martin Pinsonnault.
\newblock Packing numbers of rational ruled four-manifolds.
\newblock {\em J. Symplectic Geom.}, 11(2):269--316, 2013.

\bibitem[Bus05]{Buse05}
Olguta Buse.
\newblock Relative family {G}romov-{W}itten invariants and symplectomorphisms.
\newblock {\em Pacific J. Math.}, 218(2):315--341, 2005.

\bibitem[Bus11]{Buse11}
Olguta Buse.
\newblock Negative inflation and stability in symplectomorphism groups of ruled
  surfaces.
\newblock {\em Journal of Symplectic Geometry}, 9, 2011.

\bibitem[CCGF{\etalchar{+}}14]{CCDDHR}
Keon Choi, Daniel Cristofaro-Gardiner, David Frenkel, Michael Hutchings, and
  Vinicius Gripp~Barros Ramos.
\newblock Symplectic embeddings into four-dimensional concave toric domains.
\newblock {\em J. Topol.}, 7(4):1054--1076, 2014.

\bibitem[Che98]{Chek98}
Yu.~V. Chekanov.
\newblock Lagrangian intersections, symplectic energy, and areas of holomorphic
  curves.
\newblock {\em Duke Math. J.}, 95(1):213--226, 1998.

\bibitem[Chu12]{Church12}
Thomas Church.
\newblock Homological stability for configuration spaces of manifolds.
\newblock {\em Invent. Math.}, 188(2):465--504, 2012.

\bibitem[Cof05]{Coffey05}
Joseph Coffey.
\newblock Symplectomorphism groups and isotropic skeletons.
\newblock {\em Geom. Topol.}, 9:935--970, 2005.

\bibitem[DKK16]{DKK16}
Colin Diemer, Ludmil Katzarkov, and Gabriel Kerr.
\newblock Symplectomorphism group relations and degenerations of
  {L}andau-{G}inzburg models.
\newblock {\em J. Eur. Math. Soc. (JEMS)}, 18(10):2167--2271, 2016.

\bibitem[DRGI16]{DRGI16}
Georgios Dimitroglou~Rizell, Elizabeth Goodman, and Alexander Ivrii.
\newblock Lagrangian isotopy of tori in {$S^2\times S^2$} and {$\Bbb{C}P^2$}.
\newblock {\em Geom. Funct. Anal.}, 26(5):1297--1358, 2016.

\bibitem[Eva11a]{Eva11}
Jonathan Evans.
\newblock Symplectic mapping class groups of some stein and rational surfaces.
\newblock {\em Journal of Symplectic Geometry}, 9(1):45--82, 2011.

\bibitem[Eva11b]{Ev11}
Jonathan~David Evans.
\newblock Symplectic mapping class groups of some {S}tein and rational
  surfaces.
\newblock {\em J. Symplectic Geom.}, 9(1):45--82, 2011.

\bibitem[Gon06]{Gon06}
Eduardo Gonzalez.
\newblock Quantum cohomology and {$S^1$}-actions with isolated fixed points.
\newblock {\em Trans. Amer. Math. Soc.}, 358(7):2927--2948, 2006.

\bibitem[Gro85]{Gro85}
Misha Gromov.
\newblock Pseudoholomorphic curves in symplectic manifolds.
\newblock {\em Inventiones Mathematicae}, 82:307--347, 1985.

\bibitem[HI]{HindIvrii2}
Richard Hind and Alexander Ivrii.
\newblock Isotopies of high genus lagrangian surfaces.
\newblock https://arxiv.org/abs/math/0602475.

\bibitem[Hin04]{Hind04}
R.~Hind.
\newblock Lagrangian spheres in {$S^2\times S^2$}.
\newblock {\em Geom. Funct. Anal.}, 14(2):303--318, 2004.

\bibitem[Hin16]{Hind16}
R.~Hind.
\newblock Symplectic isotopies in dimension greater than four.
\newblock {\em J. Symplectic Geom.}, 14(4):1033--1057, 2016.

\bibitem[HK14]{HK14}
R.~Hind and E.~Kerman.
\newblock New obstructions to symplectic embeddings.
\newblock {\em Invent. Math.}, 196(2):383--452, 2014.

\bibitem[HK18]{HK18}
R.~Hind and E.~Kerman.
\newblock Correction to: {N}ew obstructions to symplectic embeddings [
  {MR}3193752].
\newblock {\em Invent. Math.}, 214(2):1023--1029, 2018.

\bibitem[HPW16]{HPW13}
Richard Hind, Martin Pinsonnault, and Weiwei Wu.
\newblock Symplectormophism groups of non-compact manifolds and space of
  lagrangians.
\newblock {\em J.Symplectic Geometry}, 1:203--226, 2016.

\bibitem[Kea14]{Ailsa14}
Ailsa~M. Keating.
\newblock Dehn twists and free subgroups of symplectic mapping class groups.
\newblock {\em J. Topol.}, 7(2):436--474, 2014.

\bibitem[Kea15]{Kea15}
Ailsa Keating.
\newblock Lagrangian tori in four-dimensional {M}ilnor fibres.
\newblock {\em Geom. Funct. Anal.}, 25(6):1822--1901, 2015.

\bibitem[KS02]{KS02}
Mikhail Khovanov and Paul Seidel.
\newblock Quivers, {F}loer cohomology, and braid group actions.
\newblock {\em J. Amer. Math. Soc.}, 15(1):203--271, 2002.

\bibitem[LL02]{LL02}
Bang-He Li and Tian-Jun Li.
\newblock Symplectic genus, minimal genus and diffeomorphisms.
\newblock {\em Asian J. Math.}, 6(1):123--144, 2002.

\bibitem[LLnt]{LL16}
Jun Li and Tian-Jun Li.
\newblock Symplectic $-2$ spheres and the symplectomorphism group of small
  rational 4-manifolds.
\newblock ArXiv Preprint.

\bibitem[LLW15]{LLW15}
Jun Li, Tian-Jun Li, and Weiwei Wu.
\newblock The symplectic mapping class group of {$\Bbb CP^2\#n\overline{\Bbb
  CP^2}$} with {$n \le 4$}.
\newblock {\em Michigan Math. J.}, 64(2):319--333, 2015.

\bibitem[LLWnt]{LLW16}
Jun Li, Tian-Jun Li, and Weiwei Wu.
\newblock Symplectic $-2$ spheres and the symplectomorphism group of small
  rational 4-manifolds, ii.
\newblock ArXiv Preprint.

\bibitem[LLWon]{LLW3}
Jun Li, Tian-Jun Li, and Weiwei Wu.
\newblock Braid groups and symplectomorphism mapping class group of rational
  surfaces.
\newblock In preparation.

\bibitem[LP04]{LP04}
Francois Lalonde and Martin Pinsonnault.
\newblock The topology of the space of symplectic balls in rational
  4-manifolds.
\newblock {\em Duke Mathematical Journal}, 122(2):347--397, 2004.

\bibitem[LW12]{LW12}
Tian-Jun Li and Weiwei Wu.
\newblock Lagrangian spheres, symplectic surface and the symplectic mapping
  class group.
\newblock {\em Geometry and Topology}, 16(2):1121--1169, 2012.

\bibitem[McD75]{McD75}
Dusa McDuff.
\newblock Configuration spaces of positive and negative particles.
\newblock {\em Topology}, 14:91--107, 1975.

\bibitem[McD98]{McD96}
Dusa McDuff.
\newblock From symplectic deformation to isotopy.
\newblock In {\em Topics in symplectic {$4$}-manifolds ({I}rvine, {CA}, 1996)},
  First Int. Press Lect. Ser., I, pages 85--99. Int. Press, Cambridge, MA,
  1998.

\bibitem[McD00]{McD00}
Dusa McDuff.
\newblock Almost complex structures on {$S^2\times S^2$}.
\newblock {\em Duke Math. J.}, 101(1):135--177, 2000.

\bibitem[McD08]{McD08}
Dusa McDuff.
\newblock The symplectomorphism group of a blow up.
\newblock {\em Geom. Dedicata}, 132:1--29, 2008.

\bibitem[MP94]{MP94}
Dusa McDuff and Leonid Polterovich.
\newblock Symplectic packings and algebraic geometry.
\newblock {\em Invent. Math.}, 115(3):405--434, 1994.
\newblock With an appendix by Yael Karshon.

\bibitem[MS98]{MS98}
Dusa McDuff and Dietmar Salamon.
\newblock {\em Introduction to symplectic topology}.
\newblock Oxford Mathematical Monographs. The Clarendon Press, Oxford
  University Press, New York, second edition, 1998.

\bibitem[MS12]{MS12}
Dusa McDuff and Felix Schlenk.
\newblock The embedding capacity of 4-dimensional symplectic ellipsoids.
\newblock {\em Ann. of Math. (2)}, 175(3):1191--1282, 2012.

\bibitem[MT06]{MT06}
Dusa McDuff and Susan Tolman.
\newblock Topological properties of {H}amiltonian circle actions.
\newblock {\em IMRP Int. Math. Res. Pap.}, pages 72826, 1--77, 2006.

\bibitem[MW]{MW15}
Cheuk~Yu Mak and Weiwei Wu.
\newblock {Dehn twists exact sequences through Lagrangian cobordism}.
\newblock {\em arXiv:1509.08028}.

\bibitem[Ops07]{Opsh07}
Emmanuel Opshtein.
\newblock Maximal symplectic packings in {$\Bbb P^2$}.
\newblock {\em Compos. Math.}, 143(6):1558--1575, 2007.

\bibitem[Ped]{Ped15}
Andres Pedroza.
\newblock Hamiltonian loops on the symplectic one-point blow up.
\newblock https://arxiv.org/abs/1510.01693v2.

\bibitem[Pin08a]{Pin08emb}
Martin Pinsonnault.
\newblock Symplectomorphism groups and embeddings of balls into rational ruled
  4-manifolds.
\newblock {\em Compos. Math.}, 144(3):787--810, 2008.

\bibitem[Pin08b]{Pin08}
Martin Pinsonnault.
\newblock Symplectomorphism groups and embeddings of balls into rational ruled
  4-manifolds.
\newblock {\em Compos. Math.}, 144(3):787--810, 2008.

\bibitem[PS08]{Pel08}
Alvaro Pelayo and Benjamin Schmidt.
\newblock Maximal ball packings of symplectic-toric manifolds.
\newblock {\em Int. Math. Res. Not. IMRN}, (3):Art. ID rnm139, 24, 2008.

\bibitem[PS16]{PS16}
Leonid Polterovich and Egor Shelukhin.
\newblock Autonomous {H}amiltonian flows, {H}ofer's geometry and persistence
  modules.
\newblock {\em Selecta Math. (N.S.)}, 22(1):227--296, 2016.

\bibitem[RB39]{RubinB}
Natanel Rubin-Blair.
\newblock The quantum johnson homomorphism and symplectomorphism of 3-folds.
\newblock http://arxiv.org/abs/1712.00339.

\bibitem[Sei97]{Seithesis}
P.~Seidel.
\newblock Floer homology and the symplectic isotopy problem.
\newblock 1997,.
\newblock Ph.D. thesis.

\bibitem[Sei99a]{Se99}
Paul Seidel.
\newblock Lagrangian two-spheres can be symplectically knotted.
\newblock {\em J. Differential Geom.}, 52(1):145--171, 1999.

\bibitem[Sei99b]{Sei99}
Paul Seidel.
\newblock On the group of symplectic automorphisms of {$\Bbb{ C} P^m \times
  \Bbb{ C} P^n$}.
\newblock In {\em Northern {C}alifornia {S}ymplectic {G}eometry {S}eminar},
  volume 196 of {\em Amer. Math. Soc. Transl. Ser. 2}, pages 237--250. Amer.
  Math. Soc., Providence, RI, 1999.

\bibitem[Sei08a]{Se4dim}
Paul Seidel.
\newblock Lectures on four-dimensional {D}ehn twists.
\newblock In {\em Symplectic 4-manifolds and algebraic surfaces}, volume 1938
  of {\em Lecture Notes in Math.}, pages 231--267. Springer, Berlin, 2008.

\bibitem[Sei08b]{Sei08}
Paul Seidel.
\newblock Lectures on four-dimensional {D}ehn twists. in.
\newblock In {\em Symplectic 4-Manifolds and Algebraic Surfaces}, pages
  231--268. volume 1938 of { Lecture {N}otes in {M}athematics}, Springer, 2008.

\bibitem[She09]{She10}
Vsevolod Shevchishin.
\newblock Secondary stiefel-whitney class and diffeomorphisms of rational and
  ruled symplectic 4-manifolds.
\newblock ArXiv preprint, 2009.

\bibitem[Sie]{Siegel16}
Kyler Siegel.
\newblock Squared dehn twists and deformed symplectic invariants.
\newblock https://arxiv.org/abs/1609.08545.

\bibitem[SS]{ShSm19}
Vsevolod Shevchishin and Gleb Smirnov.
\newblock Symplectic triangle inequality.
\newblock {\em arXiv:1908.10895}.

\bibitem[SSnt]{SS17}
Nick Sheridan and Ivan Smith.
\newblock Symplectic topology of k3 surfaces via mirror symmetry.
\newblock ArXiv Preprint.

\bibitem[Var18]{Varothesis}
Umut Varolgunes.
\newblock Mayer-vietoris property for relative symplectic cohomology.
\newblock 2018.
\newblock Ph.D. thesis.

\bibitem[Vid06]{Vid06}
Stefano Vidussi.
\newblock Lagrangian surfaces in a fixed homology class: existence of knotted
  {L}agrangian tori.
\newblock {\em J. Differential Geom.}, 74(3):507--522, 2006.

\bibitem[Wu14]{Wu14}
Weiwei Wu.
\newblock Exact {L}agrangians in {$A_n$}-surface singularities.
\newblock {\em Math. Ann.}, 359(1-2):153--168, 2014.

\bibitem[WW15]{WWfamily}
Katrin Wehrheim and Chris Woodward.
\newblock Exact triangle for fibered dehn twists.
\newblock 2015.

\bibitem[Xic02]{Xic02}
Miguel~A. Xicot\'{e}ncatl.
\newblock Product decomposition of loop spaces of configuration spaces.
\newblock In {\em Proceedings of the {F}irst {J}oint {J}apan-{M}exico {M}eeting
  in {T}opology ({M}orelia, 1999)}, volume 121, pages 33--38, 2002.

\bibitem[Yus]{Kartal18}
Baris~Kartal Yusuf.
\newblock Dynamical invariants of mapping torus categories.
\newblock https://arxiv.org/abs/1809.04046.

\bibitem[Zhaar]{ZhangJun}
Jun Zhang.
\newblock p-cyclic persistent homology and hofer distance.
\newblock {\em J. Symp. Geo.}, to appear.

\end{thebibliography}

\bibliographystyle{alpha}

\end{document}